\documentclass[12pt]{article}

\usepackage{authblk} 
\usepackage[pdfencoding=auto]{hyperref} 
\usepackage{dsfont}
\usepackage[T1]{fontenc}
\usepackage{amsmath}
\usepackage{amsfonts,amssymb}
\usepackage{extarrows}
\usepackage{geometry}
\usepackage{enumerate}
\usepackage{enumitem}

\usepackage[utf8]{inputenc}
\def\0{\emptyset}

\newtheorem{theorem}{Theorem}[section]
\newtheorem{definition}[theorem]{Definition}
\newtheorem{lemma}[theorem]{Lemma}
\newtheorem{claim}[theorem]{Claim}

\newtheorem{proposition}[theorem]{Proposition}
\newtheorem{conjecture}[theorem]{Conjecture}

\usepackage{graphicx}

\usepackage{enumerate}

\usepackage{
	amsmath,			
	amssymb,			
	enumerate,		    
	graphicx,			
	lastpage,			
	multicol,			
	multirow,			
	pifont,			    
}

\usepackage[numbers]{natbib}

\begin{document}


\title{Edge version of the inducibility  via the entropy method}

\author{\small\bf Yichen Wang\thanks{E-mail:  wangyich22@mails.tsinghua.edu.cn}}
\author{\small\bf Xiamiao Zhao\thanks{E-mail:  zxm23@mails.tsinghua.edu.cn}}
\author{\small\bf Mei Lu\thanks{E-mail:  lumei@tsinghua.edu.cn}}


\affil{\small Department of Mathematical Sciences, Tsinghua University, Beijing 100084, China}

\date{}

\maketitle

\begin{abstract}

The inducibility of a graph $H$ is about the maximum number of induced copies of $H$ in a graph on $n$ vertices.
We consider its edge version, that is, the maximum number of induced copies of $H$ in a graph with $m$ edges. Let $c(G,H)$ be the number of induced copies of $H$ in $G$ and  $\rho(H,m) = \max \{c(G,H) \mid |E(G)| = m\}$.
For any graph $H$, we prove that $\rho(H,m) = \Theta(m^{\alpha_f(H)})$ where $\alpha_f(H)$ is the fractional independence number of $H$. Therefore, we now focus on the constant factor in front of $m^{\alpha_f(H)}$.
In this paper, we give some results of $\rho(H,m)$ when $H$ is a cycle or path.
We conjecture that for any cycle $C_k$ with $k \ge 5$,  $\rho(C_k,m)= (1+o(1))\left( m/k\right)^{k/2}$ and the bound achieves by the blow up of $C_k$.
For even cycles, we establish an upper bound with an extra constant factor.
For odd cycles, we can only establish an upper bound with an extra factor depending on $k$.
We  prove that $\rho(P_{2l},m) \le \frac{m^l}{2(l-1)^{l-1}}$ and $\rho(P_{2l+1},m) \le \frac{m^{l+1}}{4l^l}$, where $l \ge 2$.
We also conjecture the asymptotic value of $\rho(P_k, m)$.
The entropy method is mainly used to prove our results.
\end{abstract}


\vskip.3cm

\section{Introduction}
Let $G$ and $H$ be two graphs.
The \textbf{induced density} of  $H$ in  $G$, denoted by $i(G,H)$, is the number of induced copies of $H$ in $G$ divided by $\binom{|V(G)|}{|V(H)|}$.
Let $i(H,n)$ denote the maximum induced density of $H$ in a graph with $n$ vertices, that is, $i(H,n) = \max \{i(G,H) \mid |V(G)| = n\}$.
The \textbf{inducibility} of  $H$, denoted by $ind(H)$, is the limit of $i(H,n)$ when $n$ goes to infinity, that is $ind(H) = \lim_{n \rightarrow +\infty}i(H,n)$.
This classic topic in extremal combinatorics is introduced by Piperger and Golumbic~\cite{pippenger1975inducibility} in 1975. In the same paper~\cite{pippenger1975inducibility}, they showed that the inducibility of every $k$-vertex  graph $H$ is at least $k!/(k^k-k)$ and conjectured that this bound is tight for a cycle $C_k$ with $k \ge 5$.

\begin{conjecture}[Pippenger and Golumbic~\cite{pippenger1975inducibility}]
The induciability of a cycle $C_k$ of length $k \ge 5$ is equal to $\frac{k!}{k^k-k}$.
\end{conjecture}

The conjecture for $k=5$ is proved by Balogh, Hu, Lidick{\`y} and Pfender~\cite{balogh2016maximum}.
For larger $k$, Piperger and Golumbic~\cite{pippenger1975inducibility} showed $ind(C_k) \le (2e+o(1))\frac{k!}{k^k}$.
The extra factor $2e$ has been improved to $128e/81$ by Hefetz and Tyomkyn~\cite{hefetz2018inducibility} in 2018.
In 2019, Král', Norin and Volec~\cite{KRAL2019359} proved that every $n$-vertex graph has at most $2n^{k}/k^k$ induced cycles of length $k$, that is, $ind(C_k) \le 2 \frac{k!}{k^k}$.

Now we fix the number of edges of graphs rather than number of vertices.
Let $c(G,H)$ be the number of induced copies of $H$ in $G$ and $\rho(H,m)$ the maximum number of induced copies of $H$ in a graph with $m$ edges. That is, $\rho(H,m) = \max \{c(G,H) \mid |E(G)| = m\}$.
If $H = K_t$, then the problem becomes the famous Kruskal-Katona Theorem~\cite{katona1987theorem, kruskal1963number} which is of central significance in extremal set theory.

Let $H$ be a graph.  Its independence number can be defined as the solution of an integer program.
Let $\omega: V(H) \rightarrow \{0,1\}$. We wish to maximize $\sum_{v \in V(H)}\omega(v)$ subject to the constraints that $\omega(u) + \omega(v) \le 1$ for all edges $uv \in E(H)$. If this integer program is relaxed such that the weights are in the interval $[0,1]$, then the solution is the \textbf{fractional independence number} of $H$, denoted by $\alpha_f(H)$. That is, maximizing $\sum_{v \in V(H)}\omega(v)$ subject to the constraints that $\omega(u) + \omega(v) \le 1$ for all edges $uv \in E(H)$ and $0 \le \omega(v) \le 1$ for all $v \in V(H)$.
Furthermore, a maximum weighting of the fractional independence number can always be obtained when restricting the weights to be in the set $\{0,1/2,1\}$~\cite{nemhauser1975vertex,willis2011bounds}. It means that equivalently we can assume $\omega(v) \in \{0,1/2,1\}$ for all $v \in V(H)$ when $\sum_{v\in V(H)} \omega(v)$ achieves maximum.
In the non-induced version of the problem. Fractional independence number plays a key role.
Alon~\cite{alon1981number}, Friegdut and Kahn~\cite{friedgut1998number} have studied the maximum number of~(not necessarily induced) copies of $G$ in a graph with $m$ edges.
\begin{theorem}[Alon~\cite{alon1981number}, Friegdut and Kahn~\cite{friedgut1998number}]\label{thm: alon friegdut kahn}
	Let $H$ be a graph without isolated vertices.
	Denote the maximum number of copies of $H$ in a graph with $m$ edges by $\pi (H, m)$. Then 
	\[
		(1-o(1)) {\left( \frac{m}{ |E(H)| }\right)}^{\alpha_f(H)} \le \pi (H,m) \le \frac{{(2m)}^{\alpha_f(H)}}{|Aut(H)|},
	\]
	where $|E(H)|$ is the number of edges in $H$, $Aut(H)$ is the automorphism group of $H$ and $|Aut(H)|$ denotes its order.
\end{theorem}

As a corollary, we have the following theorem about the induced version.

\begin{theorem}\label{thm: order}
	For any graph $H$ without isolated vertices, we have that
	\[
		(1-o(1)) {\left( \frac{m}{ |E(H)| }\right)}^{\alpha_f(H)} \le \rho (H,m) \le \frac{{(2m)}^{\alpha_f(H)}}{|Aut(H)|}.
	\]
\end{theorem}

The lower bound is achieved by the blow up of $H$ such that each vertex $v$ is blown up to ${\left(\frac{m}{|E(H)|}\right)}^{\omega(v)}$, where $\omega(v)$ is the weight of $v$ in a maximum weighting of the fractional independence number.

Note that if $H$ contains isolated vertices, the problem is meaningless since we can add infinitely many isolated vertices to $G$ such that $G$ contains infinitely many induced copies of $H$.
Theorem~\ref{thm: order} describes the asymptotic order of $\rho(H,m)$.
When $H=K_t$, $\rho(K_t,m)$ is maximized in an (almost)~complete graph $G$. At that time, $\binom{|V(G)|}{2} \approx m$, that is, $|V(G)|\approx \sqrt{2m}$, and $c(G,K_t) \approx \binom{|V(G)|}{t} \approx \binom{\sqrt{2m}}{t} \approx 2^{t/2}m^{t/2}/t!$. Note that $\alpha_f(K_t) = t/2$ and $|Aut(K_t)| = t!$. So the upper bound in Theorem~\ref{thm: order} is tight for $K_t$.

Now we focus on the constant factor in front of $m^{\alpha_f(H)}$.
We use $C_k$ to denote the cycle on $k$ vertices and $P_k$ the path on $k$ vertices.
Then  $\alpha_f(C_k) = k/2$ and  $\alpha_f(P_k) = \lfloor (k+1)/2 \rfloor$.
We have the following theorems about the number of induced copies of cycle and path.

\begin{theorem}\label{thm: odd path}We have
	$ (1+o(1))(2l+1)\frac{m^{l}}{{(2l+1)}^{l}} \le \rho(P_{2l}, m) \le \frac{m^{l}}{2{(l-1)}^{l-1}}$, where  $l \ge 2$.
\end{theorem}

\begin{theorem}\label{thm: even path} We have
	$ (1+o(1))4\frac{m^{l+1}}{{(2l+2)}^{l+1}} \le \rho(P_{2l+1}, m) \le  \frac{m^{l+1}}{4l^l}$, where $l \ge 2$.
\end{theorem}

Note that when $l=1$,  $\rho(P_3, m) \le \binom{m}{2}$ and the bound achieved by a star.

\begin{theorem}\label{thm: even cycle}
	We have the following results for even cycles.
	\begin{enumerate}
		\item We have $\rho(C_4, m) = (1+o(1))\frac{m^2}{4} $.
		\item We have $ (1+o(1)){\left( \frac{m}{6} \right)} ^{3} \le \rho(C_6, m) \le 3{\left( \frac{m}{6} \right)} ^{3}$.
		\item When $k=2l \ge 8$, we have $ (1+o(1)){\left( \frac{m}{k} \right)} ^{k/2} \le \rho(C_k, m) \le {\left( \frac{m}{2l} \right)} ^{l} {\left( 1+\frac{1}{l-1}\right)}^{l-1} \le e{\left( \frac{m}{k} \right)} ^{k/2}$.
	\end{enumerate}
\end{theorem}

\begin{theorem}\label{thm: odd cycle}
	Let $l\ge 2$. Then, we have
	\[
		(1+o(1))\left( \frac{m}{2l+1}\right)^{(2l+1)/2} \le \rho(C_{2l+1}, m) \le \frac{{(2l+1)}^{l-1/2}}{2{(l-1)}^{(l-1)(2l+1)/(2l)}} \left( \frac{m}{2l+1}\right)^{(2l+1)/2}.
	\]
\end{theorem}

We now construct the lower bound of previous theorems.
Given a graph $G$ on $s$ vertices, let $G[V_1,\ldots,V_s]$ denote the blow up of $G$ on $n\ge s$ vertices, where each vertex $v_i$ in $G$ is replaced by an independent set $V_i$ of $n_i$ vertices such that $n_1+\ldots+n_s = n$ and $u_i \in V_i$ is adjacent to $u_j \in V_j$ if and only if $v_i$ is adjacent to $v_j$ in $G$. If for any $i,j$, $|n_i - n_j| \le 1$, then we say the blow up is balanced.

\begin{figure}[t]
	\centering         
	\includegraphics[width=0.4\linewidth]{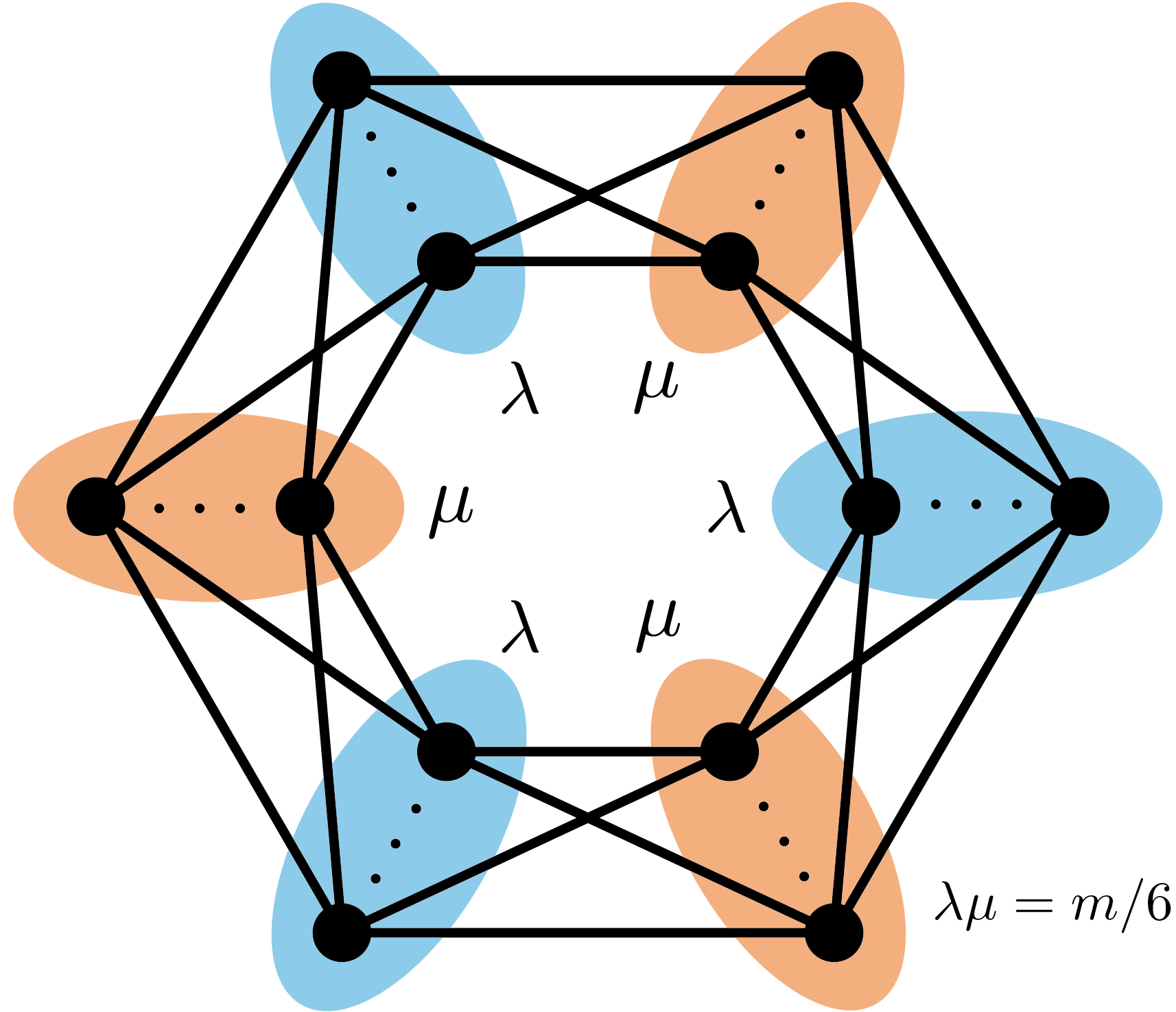}\label{1}
	\caption{Unbalanced blow up when $k=6$. The orange~(resp. blue) parts have $\mu$~(resp. $\lambda$) vertices and $\lambda \mu \approx m/6$. }\label{fig: extremal}
\end{figure}

Let $k\ge 5$. Then $\rho(C_k, m) \ge (1+o(1)) {\left( \frac{m}{k} \right)} ^{k/2}$ in Theorems \ref{thm: even cycle} and \ref{thm: odd cycle} is obtained by the construction of the balanced blow up of $C_k$, that is $C_k[V_1,\ldots,V_k]$, where each part has $n_i \approx \sqrt{m/k}$ vertices.
When $k=2l$ is even, the lower bound can also be achieved from the unbalanced blow up of $C_k$ by letting the odd parts $V_{1}, V_{3}, \ldots, V_{2l-1}$ with $\lambda$ vertices and even parts $V_{2}, V_{4}, \ldots, V_{2l}$ with $\mu$ vertices, where $\lambda \mu \approx m/k$.
When $\lambda = \mu \approx \sqrt{m/k}$, it becomes the balanced blow up of $C_k$.
Figure~\ref{1} shows an example of the unbalanced blow up of $C_6$.
We conjecture that the above lower bound is tight.

\begin{conjecture}\label{conj: cycle}
	For any $k \ge 5$,  $\rho(C_k, m) = (1+o(1)){\left( \frac{m}{k} \right)} ^{k/2}$.
\end{conjecture}

\begin{figure}[t]
	\centering         
	\includegraphics[width=0.8\linewidth]{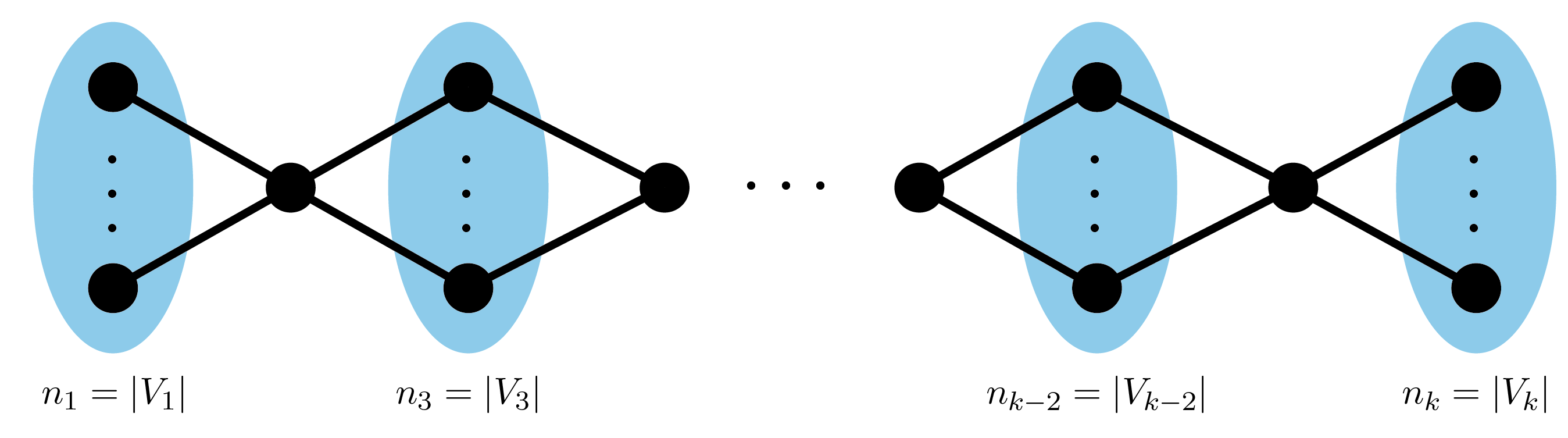}
	\caption{The blow up of $P_k$ when $k$ is odd.}\label{fig: even path}
\end{figure}

Now we give the lower bounds of $\rho(P_{k}, m)$ in Theorems \ref{thm: odd path} and \ref{thm: even path}. Note that $\rho(P_{k}, m) \ge (1+o(1))(k+1) \left( \frac{m}{k+1} \right)^{\frac{k}{2}}$ where $k$ is an even number and the bound achieved by the balanced blow up of $C_{k+1}$.
When $k$ is an odd number, we have that $\rho(P_{k}, m) \ge (1+o(1))4 {\left( \frac{m}{k+1} \right)}^{(k+1)/2}$ by the construction of the unbalanced blow up of $P_k$~(see Figure~\ref{fig: even path}).
Let $P_k[V_1,\ldots,V_k]$ be a blow up of $P_k$ and denote $n_i = |V_{i}|, 1 \le i \le k$.
Let $n_i = 1$ for even $i$. Then $m = n_1 + 2n_3 + 2n_5 + \ldots + 2n_{k-2} + n_{k}$ and the number of induced copies of $P_k$ is $n_1 n_3\ldots n_k \le 4 {(m/(k+1))}^{(k+1)/2}$ and the equality holds when $n_1 = 2n_3 = 2n_5 = \ldots = 2n_{k-2} = n_{k}$. So we conjecture those bounds to be tight.

\begin{conjecture}\label{conj: odd path}
	Let $k \ge 4$ be an even integer, then $\rho(P_{k}, m) = (1+o(1)) (k+1) \left( \frac{m}{k+1} \right)^{\frac{k}{2}}$.
\end{conjecture}

\begin{conjecture}\label{conj: even path}
	Let $k \ge 5$ be an odd integer, then $\rho(P_{k}, m) = (1+o(1))4 {\left( \frac{m}{k+1} \right)}^{(k+1)/2}$.
\end{conjecture}

Comparing Theorem~\ref{thm: odd path} with Conjecture~\ref{conj: odd path},  our theorem have an extra factor $\frac{{(2l+1)}^{l-1}}{2{(l-1)}^{l-1}} = \Theta(2^{l-2}e^{3/2})$ when $l$ goes to infinity. Similarly, comparing Theorem~\ref{thm: odd cycle} to Conjecture~\ref{conj: cycle}, there is an extra factor $\frac{{(2l+1)}^{l-1/2}}{2{(l-1)}^{(l-1)(2l+1)/(2l)}} = \Theta(2^{l-3/2}e^{3/2})$ when $l$ goes to infinity.
That is because Theorem~\ref{thm: odd cycle} is actually a corollary from Theorem~\ref{thm: odd path} using Shearer's Lemma~(Lemma~\ref{lem: shearer}) which we will mention in Section~\ref{sec: preliminary}. The extra error factor when estimating induced odd cycle is due to our overestimation about induced odd path.

We use the entropy method to prove our main results.
The entropy method is a powerful tool in such problems, for example, counting rainbow triangles when fixing the number of edges in each color class~\cite{chao2024kruskal}; counting $t$-cliques when fixing the $p$-norm of the degree sequence~\cite{chao2024many}.

This paper is organized as follows. In Section~\ref{sec: preliminary},  we introduce basic properties of the entropy method.
In Section~\ref{sec: even path}, we will prove Theorem~\ref{thm: even path}.
In Section~\ref{sec: odd cycle odd path}, we will prove Theorems~\ref{thm: odd path} and \ref{thm: odd cycle}.
In Section~\ref{sec: even cycle}, we will prove Theorem \ref{thm: even cycle}.

\section{Preliminary}\label{sec: preliminary}

In this section, we introduce basic properties of the entropy method.
For a complete discussion and proofs for properties of the entropy method, readers may refer to Section 14.6 in~\cite{alon2016probabilistic}.

\begin{definition}\label{def: entropy}
	For any discrete random varible $X$ with finite support $S$, the \textbf{entropy} of $X$ is defined as
	\begin{equation*}
			H(X) = -\sum_{x \in S} p_X(x) \log p_X(x),
	\end{equation*}
	where $p_X(x)$ denotes  $\mathbb{P}(X=x)$.
\end{definition}

The entropy $H(X_1, \ldots, X_n)$ of several discrete random variables $X_1, \ldots, X_n$ is defined as the entropy of the random variable $X=(X_1, \ldots, X_n)$.

\begin{definition}[Conditional Entropy]\label{def: conditional entropy}
	Let $X,Y$ be two random variables with finite supports $S,T$, respectively.
	Denote by $p_{X,Y}(x,y)$ the probability $\mathbb{P}(X=x,Y=y)$.
	For each $y \in T$, let $X|Y=y$ be the random variable $X$ conditioned on $Y=y$, and $H(X \mid Y=y)$ be its entropy.
	The conditional entropy of $X$ given $Y$ is defined as
	\begin{equation*}
		H(X \mid Y) = \sum_{y \in T} p_Y(y) H(X\mid Y=y) = -\sum_{x \in S, y \in T} p_{X,Y}(x,y) \log \left(\frac{p_{X,Y}(x,y)}{p_Y(y)}\right).
	\end{equation*}
\end{definition}

\begin{proposition}\label{prop: entropy support bound}
	For any random variable $X$ with finite support $S$, we have $H(X) \le \log |S|$.
	The equality holds when $X$ is uniformly distributed over $S$.
\end{proposition}

\begin{proposition}[Chain Rule]\label{prop: chain rule}
	For any two random variables $X$ and $Y$ with finite supports, we have
	\begin{equation*}
		H(X,Y) = H(X) + H(Y\mid X) = H(Y) + H(X\mid Y).
	\end{equation*}
	More generally, for any $n$ random variables $X_1, \ldots, X_n$ with finite supports, we have
	\begin{equation*}
		H(X_1, \ldots, X_n) = H(X_1) + H(X_2\mid X_1) + \cdots + H(X_{n} \mid X_{1}, \ldots, X_{n-1}).
	\end{equation*}
\end{proposition}

\begin{proposition}[Drop Condition]\label{prop: drop condition}
	For any two random variables $X$ and $Y$ with finite supports, we have
	\[H(X \mid Y) \le H(X).\]
\end{proposition}

\begin{lemma}[Shearer's Lemma]\label{lem: shearer}
	Let $X_1,  \ldots, X_n$ be random variables with finite supports. Let $\mathcal{A}$ be a collection of subsets of $[n]$ such that each $i \in [n]$ is contained in at least $r$ sets in $\mathcal{A}$. Then
	\[
	H(X_1, X_2, \ldots, X_n) \le \frac{1}{r} \sum_{A \in \mathcal{A}} H((X_i)_{i \in A}).
	\]
\end{lemma}

We say an ordered edge tuple $(e_1,e_2,\ldots,e_k)$ is an ordered induced copy of $C_k$ if $e_i$ intersects $e_{i-1}$ and $e_{i+1}$ on two different vertices for $i=1,\ldots, k$~(let $e_0 = e_k, e_{k+1} = e_1$) and there is no chord inside.
Similarly, we say $(e_1,e_2,\ldots,e_{k-1})$ is an ordered induced copy of $P_k$ if $e_i$ intersects $e_{i-1}$ and $e_{i+1}$ on two different vertices for $i=2,\ldots, k-2$ and there is no other edges connecting $e_i$ and $e_j$ for any $1 \le i \neq j \le k-1$.

Given a graph $G$, we say a sequence of edges $(e_1,e_2,\ldots,e_t)$ is \textbf{well-ordered} if $v_iu_{i+1}\in E(G)$  for all $1 \le i \le t-1$ but there is no other edge among $\{u_iu_j, u_iv_j, v_iu_j, v_iv_j\}_{i \neq j}$, where  $e_i = u_iv_i, 1 \le i \le t$.
For example, in an ordered induced copy of $C_k = (e_1,\ldots,e_k)$, $(e_1,e_3,\ldots,e_{2t-1})$~($2t-1 < k-1$) is well-ordered.
In the following, when there is no ambiguity, we usually ignore specifying the graph $G$.

We say a well-ordered edge tuple $(e_1,e_3,\ldots,e_{2t-1})$ is \textbf{embedded} in an ordered induced copy of $P_k = (e_1',e_2',\ldots, e_{k-1}')$ if $e_{2i-1} = e_{2i-1}'$ for all $1 \le i \le t$. For convenience, we write $e^{o}_{<2t+1} = e_1, e_3, \ldots, e_{2t-1}$. Let $e^{o}_{<2t+1} = e_1, e_3, \ldots, e_{2t-1}$ and $\tilde{e}^{o}_{<2t+1} = \tilde{e}_1, \tilde{e}_3, \ldots, \tilde{e}_{2t-1}$.
We say $e^{o}_{<2t+1} = \tilde{e}^{o}_{<2t+1}$ if $e_i = \tilde{e}_i$ for $i=1,3,\ldots, 2t-1$.

\section{Proof of Theorem~\ref{thm: even path}}\label{sec: even path}

Let $k = 2l+1 \ge 5$ be an odd integer and $G$ the extremal graph maximizing the number of induced copies of $P_k$.
Denote the number of induced copies of $P_k$ in $G$ by $\Upsilon$.
Let $P = (e_1,e_2,\ldots, e_{k-1})$ be an random variable uniformly chosen from all ordered induced copies of $P_k$ in $G$. Hence, we have $H(P) = \log (2\Upsilon)$ by Proposition~\ref{prop: entropy support bound}.
For a well-ordered edge tuple $(e_1,e_3,\ldots,e_{2t-1})$~($1 \le t \le l-2)$, let
$\beta(e_1,e_3,\ldots,e_{2t-1})$ be the number of ordered induced copies of $P_{k} = (e_1',e_2',\ldots,e_k')$ such that $e_{2i-1}=e_{2i-1}'$ for all $1 \le i \le t$ and $\alpha(e_1,e_3,\ldots,e_{2t-1})$ be the number of edges $e$ such that $(e_1,e_3,\ldots,e_{2t-1},e)$ is well-ordered.
For a well-ordered edge tuple $(e_1,e_3,\ldots,e_{2l-3})$, let $\mathcal{S}(e_1,e_3,\ldots,e_{2l-3})$ be an edge set contained $e$ such that there exists an ordered induced path $P_k = (e_1',e_2',\ldots,e_{2l}')$ satisfying $e_{2i-1}=e_{2i-1}'$ for all $1 \le i \le l-1$ and $e_{2l}'=e$.
Write $\gamma_{0}(e_1,e_3,\ldots,e_{2l-3}) = |\mathcal{S}(e_1,e_3,\ldots,e_{2l-3})|$.
For a well-ordered edge tuple $(e_1,e_3,\ldots,e_{2l-3})$ and an edge $e_{2l} \in \mathcal{S}(e_1,e_3,\ldots,e_{2l-3})$, let $\gamma_1(e_1,e_3,\ldots,e_{2l-3}, e_{2l})$~(resp. $\gamma_2(e_1,e_3,\ldots,e_{2l-3}, e_{2l})$) be the number of $e$ such that there exists an ordered induced path $P_k = (e_1',e_2',\ldots,e_{2l}')$ satisfying $e_{2i-1}=e_{2i-1}'$ for all $1 \le i \le l-1$, $e_{2l}'=e_{2l}$ and $e_{2l-2}'=e$~(resp. $e_{2l-1}' = e$).

Note that $H(e_2,e_4,\ldots,e_{2l-4} \mid e^{o}_{<2l-1}) = 0$ by Proposition~\ref{prop: entropy support bound} since $e_2,e_4,\ldots,e_{2l-4}$ are uniquely determined by $e^{o}_{<2l-1}$. Moreover, we have $H(e_{2l-2},e_{2l-1} \mid e^{o}_{<2l-1}, e_{2l}) = H(e_{2l-2} \mid e^{o}_{<2l-1}, e_{2l}) = H(e_{2l-1} \mid e^{o}_{<2l-1}, e_{2l})$ since $e_{2l-2}$~(resp. $e_{2l-1}$) is uniquely determined by $e^{o}_{<2l-1}, e_{2l}$ and $e_{2l-1}$~(resp. $e_{2l-2}$).
Combining Propositions~\ref{prop: chain rule} and \ref{prop: entropy support bound}, we have

\begin{equation}\label{eq: even path H(P)}
\begin{aligned}
	&H(P) = \log (2\Upsilon) \\
	\le & H(e_1) + \sum_{i=1}^{l-2} H(e_{2i+1} \mid e^{o}_{<2i+1}) + H(e_{2l} \mid e^{o}_{<2l-1})  \\
	 & + \frac{1}{2} \left( H(e_{2l-2} \mid e^{o}_{<2l-1}, e_{2l}) + H(e_{2l-1} \mid e^{o}_{<2l-1}, e_{2l}) \right).
\end{aligned}
\end{equation}

Now we analyze each term in (\ref{eq: even path H(P)}). Firstly, we have

\begin{equation} \label{eq: even path 2i+1 term}
\begin{aligned}
	H(e_{2i+1} \mid e^{o}_{<2i+1}) & \le \sum_{\tilde{e}^{o}_{<2i+1}} \mathbb{P}(e^{o}_{<2i+1} = \tilde{e}^{o}_{<2i+1}) \log \alpha(\tilde{e}^o_{<2i+1}) \\
	& = \sum_{\tilde{e}^{o}_{<2i+1}} \frac{\beta(\tilde{e}^o_{<2i+1})}{2\Upsilon} \log \alpha(\tilde{e}^o_{<2i+1}) \\
	& = \sum_{\tilde{e}^{o}_{<2i+1}} \sum_{\text{Ordered induced $P_k$: $\tilde{P}$}} \frac{\mathrm{1}_{\text{ $\tilde{e}^o_{<2i+1}$ is embedded in $\tilde{P}$}}}{2\Upsilon} \log \alpha(\tilde{e}^o_{<2i+1}) \\
	& = \sum_{\text{Ordered induced $P_k=(\tilde{e}_1,\ldots,\tilde{e}_{k-1})$}} \frac{1}{2\Upsilon} \log \alpha(\tilde{e}^o_{<2i+1}).
\end{aligned}
\end{equation}

Similarly, we have
\begin{equation*}
\begin{aligned}
	H(e_{2l} \mid e^{o}_{<2l-1}) &\le \sum_{\text{Ordered induced $P_k=(\tilde{e}_1,\ldots,\tilde{e}_{k-1})$}} \frac{1}{2\Upsilon} \log \gamma_0(\tilde{e}^o_{<2l-1}), \\
	H(e_{2l-2} \mid e^{o}_{<2l-1}, e_{2l}) &\le \sum_{\text{Ordered induced $P_k=(\tilde{e}_1,\ldots,\tilde{e}_{k-1})$}} \frac{1}{2\Upsilon} \log \gamma_1(\tilde{e}^o_{<2l-1}, \tilde{e}_{2l}), \\
	H(e_{2l-1} \mid e^{o}_{<2l-1}, e_{2l}) &\le \sum_{\text{Ordered induced $P_k=(\tilde{e}_1,\ldots,\tilde{e}_{k-1})$}} \frac{1}{2\Upsilon} \log \gamma_2(\tilde{e}^o_{<2l-1}, \tilde{e}_{2l}).
\end{aligned}
\end{equation*}

When fixing the ordered induced path $P_k=(\tilde{e}_1,\ldots,\tilde{e}_{k-1})$, $\{\alpha(\tilde{e}^o_{<2i+1})\}_{i=1,2,\ldots,l-2}$, $\gamma_0(\tilde{e}^o_{<2l-1})$, $\gamma_1(\tilde{e}^o_{<2l-1}, \tilde{e}_{2l})$ and $\gamma_2(\tilde{e}^o_{<2l-1}, \tilde{e}_{2l})$ are all cardinality of some edge sets. By the definitions, those edge sets are disjoint. Hence, we have
\[
\sum_{i=1}^{l-2}\alpha(\tilde{e}^o_{<2i+1}) + \gamma_0(\tilde{e}^o_{<2l-1}) + \gamma_1(\tilde{e}^o_{<2l-1}, \tilde{e}_{2l}) + \gamma_2(\tilde{e}^o_{<2l-1}, \tilde{e}_{2l}) \le m.
\]

By AM-GM inequality, we have that
\begin{equation*}
\begin{aligned}
	& \prod_{i=1}^{l-2} \alpha(\tilde{e}^o_{<2i+1})^2 \cdot \gamma_0(\tilde{e}^o_{<2l-1})^2 \cdot {\gamma_1(\tilde{e}^o_{<2l-1}, \tilde{e}_{2l})} \cdot {\gamma_2(\tilde{e}^o_{<2l-1}, \tilde{e}_{2l})} \\
	\le & \frac{1}{4}  {\left( \frac{2\left( \sum_{i=1}^{l-2}\alpha(\tilde{e}^o_{<2i+1}) + \gamma_0(\tilde{e}^o_{<2l-1}) + \gamma_1(\tilde{e}^o_{<2l-1}, \tilde{e}_{2l}) + \gamma_2(\tilde{e}^o_{<2l-1}, \tilde{e}_{2l}) \right)}{2l} \right)}^{2l} \\
	\le & \frac{1}{4}  {\left( \frac{2m}{2l} \right)}^{2l} = \frac{1}{4}  {\left( \frac{m}{l} \right)}^{2l}. \\
\end{aligned}
\end{equation*}

Put above inequalities into~\eqref{eq: even path H(P)}, we have

\begin{equation*}
\begin{aligned}
	&H(P) = \log (2\Upsilon) \\
	\le & \log m + \sum_{\text{Ordered induced $P_k=(\tilde{e}_1,\ldots,\tilde{e}_{k-1})$}} \frac{1}{4\Upsilon} \log \left( \frac{1}{4}  {\left( \frac{m}{l} \right)}^{2l} \right) \\
	\le & \log \frac{m^{l+1}}{2 l^l},
\end{aligned}
\end{equation*}
which implies Theorem~\ref{thm: even path} holds.

\section{Proofs of Theorems~\ref{thm: odd path} and \ref{thm: odd cycle}}\label{sec: odd cycle odd path}

We first prove Theorem~\ref{thm: odd cycle} when assuming Theorem~\ref{thm: odd path} holds.

\vspace{.2cm}
\noindent\textbf{Proof of Theorem~\ref{thm: odd cycle}:}

Let $G$ be the extremal graph with $m$ edges maximizing the number of induced copies of $C_{2l+1}$, where $l\ge 2$.
We say an ordered vertex tuple $(v_1,v_2,\ldots,v_{2l+1})$ is an induced copy of  $C_{2l+1}$ if $v_iv_{i+1}\in E(G)$ for all $1 \le i \le 2l+1$~(let $v_{2l+2} = v_1$), and there are no other edges among these vertices.
In other words, $v_1v_2\ldots v_{2l+1}v_1$ is an induced cycle in $G$.
Similarly, we say an ordered vertex tuple $(v_1,v_2,\ldots,v_{2l})$ is an induced copy of $P_{2l}$ if $v_iv_{i+1}\in E(G)$ for all $1 \le i \le 2l-1$, and there are no other edges among these vertices.
Let $C=(v_1,v_2,\ldots,v_{2l+1})$ be a random vertex tuple uniformly chosen from all ordered induced copies of $C_{2l+1}$ in $G$.
Then the entropy of $C$ is $H(C) = \log (2(2l+1)\Gamma )$, where $\Gamma$ denotes the number of induced copies of $C_{2l+1}$ in $G$.

By Shearer's Lemma~(Lemma~\ref{lem: shearer}), we have
\begin{equation*}
\begin{aligned}
	H(C) & \le \frac{1}{2l} \sum_{i=1}^{2l+1}H(v_1, \ldots, v_{i-1}, v_{i+1},\ldots, v_{2l+1}). \\
\end{aligned}
\end{equation*}

Note that $(v_{i+1}, \ldots, v_{2l+1}, v_{1}, \ldots, v_{i-1})$ is always an ordered induced copy of $P_{2l}$ in $G$.
Let $\Upsilon$ denote the number of induced copies of $P_{2l}$ in $G$.
Then the number of ordered $P_{2l}$ is $2\Upsilon$ and we have that $H(v_{i+1}, \ldots, v_{2l+1}, v_{1}, \ldots, v_{i-1}) \le \log (2\Upsilon)$ for any $1 \le i \le 2l+1$.
Consequently, we have
\begin{equation*}
	H(C) = \log (2(2l+1)\Gamma )  \le \frac{2l+1}{2l} \log (2\Upsilon).
\end{equation*}

It is worth mentioning that the above inequality is not tight in our conjectured extremal graph, blow up of $C_{2l+1}$.
Combining Theorem~\ref{thm: odd path}, we have that
\begin{equation*}
\begin{aligned}
	\Gamma & \le \frac{1}{2(2l+1)} {(2\Upsilon)}^{\frac{2l+1}{2l}} \\
	& \le \frac{{(2l+1)}^{l-1/2}}{2{(l-1)}^{(l-1)(2l+1)/(2l)}} \left( \frac{m}{2l+1} \right)^{(2l+1)/2}, \\
\end{aligned}
\end{equation*}
which proves Theorem~\ref{thm: odd cycle}. \hfill $\square$

In the following, we give the proof of Theorem~\ref{thm: odd path}.

\vspace{.2cm}
\noindent\textbf{Proof of Theorem~\ref{thm: odd path}:}
Let $G$ be the extremal graph with $m$ edges maximizing the number of induced copies of $P_{2l}$. Let $\Upsilon$ denote the number of induced copies of $P_{2l}$ in $G$.
For a well-ordered edge tuple $(e_1,e_3,\ldots,e_{2t-1})$~($1 \le t \le l-2)$, let $\beta(e_1,e_3,\ldots,e_{2t-1})$ be the number of ordered induced copies of $P_{2l} = (e_1',e_2',\ldots,e_{2l}')$ such that $e_{2i-1}=e_{2i-1}'$ for all $1 \le i \le t$ and $\alpha(e_1',e_3',\ldots,e_{2t-1}')$ be the number of edges $e$ such that $e_1',e_3',\ldots,e_{2t-1}', e$ is well-ordered.

Let $P = (e_1, e_2,\ldots,e_{2l-1})$ be a random variable uniformly chosen from all ordered induced copies of $P_{2l}$ in $G$.
Then the entropy of $P$ is $H(P) = \log (2\Upsilon)$, and from Proposition~\ref{prop: chain rule} we have
\begin{equation*}
\begin{aligned}
	H(P) & = H(e_1) + \sum_{i=1}^{l-1}H(e_{2i+1} \mid e^{o}_{<2i+1}) + H(e_2,e_4,\ldots,e_{2l-2} \mid e_1,e_3,\ldots,e_{2l-1}). \\
\end{aligned}
\end{equation*}

Firstly, note that $H(e_2,\ldots,e_{2l-2} \mid e_1,e_3,\ldots,e_{2l-1}) = 0$ by Proposition~\ref{prop: entropy support bound} since $e_2,\ldots,e_{2l-2}$ are uniquely determined by $e_1,e_3,\ldots,e_{2l-1}$.
Moreover, similar to Eq~\ref{eq: even path 2i+1 term}, we have

\begin{equation*}
\begin{aligned}
	H(e_{2i+1}\mid e^{o}_{< 2i+1}) \le \sum_{\text{Ordered induced $\tilde{P}_{2l}=\tilde{e}_1 \tilde{e}_2\ldots \tilde{e}_{2l-1}$}} \frac{1}{2\Upsilon} \left(\log \alpha(\tilde{e}_{1}, \tilde{e}_{3}, \ldots, \tilde{e}_{2i-1}) \right).\\
\end{aligned}
\end{equation*}

For convenience, write $\alpha_i = \alpha(\tilde{e}_{1}, \tilde{e}_{3}, \ldots, \tilde{e}_{2i-1})$, and then we have
\begin{equation*}
\begin{aligned}
	&H(P) = H(e_1) + \sum_{i=1}^{l-1}H(e_{2i+1} \mid e^{o}_{<2i+1}) \\
	\le&  \log m + \sum_{\text{Ordered induced~} P_{2l}=\tilde{e}_1 \tilde{e}_2\ldots \tilde{e}_{2l-1}} \sum_{i=1}^{l-1} \frac{1}{2\Upsilon} \left(\log \alpha_i \right) \\
	\le & \log m + \sum_{\text{Ordered induced~} P_{2l}=\tilde{e}_1 \tilde{e}_2\ldots \tilde{e}_{2l-1}} \frac{1}{2\Upsilon} \left(\log \prod_{i=1}^{l-1} \alpha_i \right) \\
	\le & \log m + \sum_{\text{Ordered induced~} P_{2l}=\tilde{e}_1 \tilde{e}_2\ldots \tilde{e}_{2l-1}} \frac{1}{2\Upsilon} \left(\log \left( \frac{\sum_{i=1}^{l-1} \alpha_i }{l-1} \right)^{l-1} \right). \\
\end{aligned}
\end{equation*}

Note that each $\alpha_i$ is a cardinality of some edge set, and those edge sets are disjoint from each other by their definitions. Therefore, we have that $\sum_{i=1}^{l-1} \alpha_i \le m$.
Thus, we obtain

\begin{equation*}
\begin{aligned}
	H(P) = \log (2\Upsilon) & \le \log m + \sum_{\text{Ordered induced~}P_{2l}=\tilde{e}_1 \tilde{e}_2\ldots \tilde{e}_{2l}} \frac{1}{2\Upsilon} \log \left( \frac{m}{l-1} \right)^{l-1} \\
	& = \log m + \log \left( \frac{m}{l-1} \right)^{l-1} \\
	& = \log \frac{m^{l}}{{(l-1)}^{l-1}},\\
\end{aligned}
\end{equation*}
which proves our theorem. \hfill $\square$

\section{Proof of Theorem~\ref{thm: even cycle}}\label{sec: even cycle}

Assume $k=2l$, we say an ordered edge tuple $(e_1, e_2, \ldots, e_l)$ \textbf{characterizes} an induced copy of $C_k$ if we write $e_i = u_iv_i$ and then $v_iu_{i+1}$ is an edge for $i=1,\ldots, l$~(let $u_{l+1} = u_0$) and there is no other edge among $\{u_iu_j, u_iv_j, v_iu_j, v_iv_j\}_{i \neq j}$. For example, if $(e_1,e_2,\ldots,e_{2l})$ is an ordered induced copy of $C_{2l}$, then $(e_1,e_3,e_5,\ldots,e_{2l-1})$ characterizes an ordered induced $C_{2l}$.



For a well-ordered edge tuple $(e_1,e_3,\ldots,e_{2t-1})$~($1 \le t \le l-1$), let $\beta(e_1,\ldots,e_{2t-1})$ be the number of ordered induced copies of $C_{k} = (e_1',e_2',\ldots,e_k')$ such that $e_{2i-1}=e_{2i-1}'$ for all $1 \le i \le t$.
When $t \le l-2$, let $\alpha(e_1,e_3,\ldots,e_{2t-1})$ be the number of edges $e$ such that $(e_1,e_3,\ldots,e_{2t-1},e)$ is well-ordered. When $t= l-1$, let $\alpha(e_1,e_3,\ldots,e_{2l-3})$ be the number of edges $e$ such that $(e_1,e_3,\ldots,e_{2l-3},e)$ characterizes an induced copy of $C_k$.


\subsection{Case 1: \texorpdfstring{$k=4$}{k=4}}

We first prove the upper bound $\rho(C_4,m) \le \frac{m^2}{4}$.
For any graph $G$, let $C = (e_1,e_2,e_3,e_4)$ be a random variable uniformly distributed among all ordered induced copies of $C_4$.
Denote the number of induced copies of $C_4$ in $G$ by $\Gamma$.
Then $H(C) = \log (8\Gamma) \le H(e_1) + H(e_3) + H(e_2,e_4 \mid e_1,e_3)\le 2\log m + \log 2$ and hence, $\Gamma \le \frac{m^2}{4}$.
Here, $H(e_1), H(e_3) \le \log m$ is from Proposition~\ref{prop: entropy support bound} and $H(e_2,e_4 \mid e_1,e_3) = 2$ is because once $e_1$ and $e_3$ are determined, there are only two possible choices for $(e_2,e_4)$.

The lower bound $\rho(C_4,m) \ge (1+o(1))\frac{m^2}{4}$ is achieved by a balanced complete bipartite graph $H$, where each part have $n \approx \sqrt{m}$ vertices. Then the number of induced copies of $C_4$ in $H$ is $\binom{n}{2}^2 = (1+o(1))\frac{m^2}{4}$.

\subsection{Case 2: \texorpdfstring{$k\ge 8$}{k>=8}}\label{subsec: main}

In this subsection, assume $k=2l \ge 8$.
Let $G$ be the extremal graph maximizing the number of induced copies of $C_{k}$.
Let $C=(e_1,\ldots,e_{2l})$ be a random variable uniformly chosen from all ordered induced copies of $C_{k}$ in $G$.
Then the entropy of $C$ is $H(C) = \log (2k\Gamma ) = \log (4l\Gamma )$, where $\Gamma$ denotes the number of induced copies of $C_{2l}$ in $G$.

By Proposition~\ref{prop: chain rule}, we have
\begin{equation*}
\begin{aligned}
	H(e_1,\ldots, e_{2l}) & \le H(e_1) + H(e_{3} \mid e_1) + \sum_{i=2}^{l-1} H(e_{2i+1}\mid e^{o}_{< 2i+1}) + H(e_2,\ldots, e_{2l} \mid e^{o}_{< 2l+1}). \\
\end{aligned}
\end{equation*}

Note that once we determine $e_1, e_3, \ldots, e_{2l-1}$, the remaining edges $e_2, e_4, \ldots, e_{2l}$ are uniquely determined since $C$ is an induced cycle. Hence, by Proposition~\ref{prop: entropy support bound}, we have $H(e_2,\ldots, e_{2l} \mid e^{o}_{< 2l+1}) = 0$.
Since $(e_1,e_3)$ is well-ordered in $C$, we have
\begin{equation*}
	\begin{aligned}
		H(e_3 \mid e_1) & \le \sum_{\tilde{e}_1} \mathbb{P}(e_1=\tilde{e}_1) \log \alpha(\tilde{e}_1)  \\
		& = \sum_{\tilde{e}_1} \frac{\beta(\tilde{e}_1)}{4l\Gamma} \log \alpha(\tilde{e}_1)  \\
		& =  \sum_{\tilde{e}_1}  \sum_{\tilde{C} \text{ is an induced $C_{2l}$ } } \frac{\mathrm{1}_{\text{$\tilde{e}_1$ is in $\tilde{C}$}}}{2l\Gamma} \log \alpha(\tilde{e}_1)  \\
		& =  \sum_{C_{2l}=\tilde{e}_1 \tilde{e}_2\tilde{e}_3\ldots \tilde{e}_{2l}\tilde{e}_1}  \frac{1}{2l\Gamma} \sum_{i=1}^{2l} \log \alpha(\tilde{e}_i).  \\
	\end{aligned}
\end{equation*}
Here the second equality is due to for each induced cycle $\tilde{C}$ contained  $\tilde{e}_1$, there are two ordered $\tilde{C}$ such that $\tilde{e}_1$ is the first coordinate.
Similarly, note that $e^{o}_{< 2i+1}$ is well-ordered in $C$ for $i=2,\ldots, l-1$. Then we have
\begin{equation*}
	\begin{aligned}
		&H(e_{2i+1}\mid e^{o}_{< 2i+1}) \\
		\le&  \sum_{\tilde{e}^{o}_{< 2i+1}}  \mathbb{P}(e^{o}_{< 2i+1}=\tilde{e}^{o}_{< 2i+1}) \log \alpha(\tilde{e}^{o}_{< 2i+1}) \\
		=&  \sum_{\tilde{e}^{o}_{< 2i+1}} \frac{\beta(\tilde{e}^{o}_{< 2i+1})}{4l\Gamma} \log \alpha(\tilde{e}^{o}_{< 2i+1}) \\
		=&  \sum_{C_{2l}=\tilde{e}_1 \tilde{e}_2\ldots \tilde{e}_{2l}\tilde{e}_1} \sum_{\tilde{e}^{o}_{< 2i+1}} \frac{1}{4l\Gamma} \log \alpha(\tilde{e}^{o}_{< 2i+1}) \\
		=&  \sum_{C_{2l}=\tilde{e}_1 \tilde{e}_2\ldots \tilde{e}_{2l}\tilde{e}_1} \frac{1}{4l\Gamma} \left(\sum_{j=1}^{2l} \log \alpha(\tilde{e}_{j}, \tilde{e}_{j+2}, \ldots, \tilde{e}_{j+2(i-1)})+\sum_{j=1}^{2l} \log \alpha(\tilde{e}_{j}, \tilde{e}_{j-2}, \ldots, \tilde{e}_{j-2(i-1)}) \right).\\
	\end{aligned}
\end{equation*}

For convenience, write $\alpha^{+}_{j,i} = \alpha(\tilde{e}_{j}, \tilde{e}_{j+2}, \ldots, \tilde{e}_{j+2(i-1)})$ and $\alpha^{-}_{j,i} = \alpha( \tilde{e}_{j}, \tilde{e}_{j-2}, \ldots, \tilde{e}_{j-2(i-1)})$ for $2 \le i \le l-1$ and $1 \le j \le 2l$.
Here both $\alpha^{+}_{j,i}$ and $\alpha^{-}_{j,i}$ depend on some induced cycle copy of $C_{2l} =\tilde{e}_1\tilde{e}_2\ldots \tilde{e}_{2l}\tilde{e}_1$ which is not reflected in the notation.
Hence in the summation below, the value of $\alpha^{+}_{j,i}$ and $\alpha^{-}_{j,i}$ may change under different induced cycle copies of $C_{2l}$.
Combining those inequalities above, we have
\begin{equation*}
\begin{aligned}
	&\log (4l\Gamma) = H(e_1, \ldots, e_{2l})\\
	\le&  \log m +  \sum_{C_{2l}=\tilde{e}_1 \tilde{e}_2\ldots \tilde{e}_{2l}\tilde{e}_1}  \frac{1}{2l\Gamma} \sum_{j=1}^{2l} \log \alpha(\tilde{e}_j) + \sum_{i=2}^{l-1}\sum_{C_{2l}=\tilde{e}_1 \tilde{e}_2\ldots \tilde{e}_{2l}\tilde{e}_1} \frac{1}{4l\Gamma} \left(\sum_{j=1}^{2l} \log \alpha^{+}_{j,i}+\sum_{j=1}^{2l} \log \alpha^{-}_{j,i} \right) \\
	\le& \log m + \sum_{C_{2l}=\tilde{e}_1 \tilde{e}_2\ldots \tilde{e}_{2l}\tilde{e}_1}\frac{1}{4l\Gamma}\left( \sum_{j=1}^{2l} \left( \log \alpha(\tilde{e}_j) +  \sum_{i=2}^{l-1}  \log \alpha^{+}_{j,i} \right)+\left( \log \alpha(\tilde{e}_j) +  \sum_{i=2}^{l-1}  \log \alpha^{-}_{j,i} \right) \right) \\
	\le&  \log m + \sum_{C_{2l}=\tilde{e}_1 \tilde{e}_2\ldots \tilde{e}_{2l}\tilde{e}_1}\frac{1}{4l\Gamma}\left(  \log \left( \prod_{j=1}^{2l} \left(2\frac{\alpha(\tilde{e}_j)}{2} \prod_{i=2}^{l-1}  \alpha^{+}_{j,i}\right) \right)+\log \left( \prod_{j=1}^{2l} \left(2\frac{\alpha(\tilde{e}_j)}{2} \prod_{i=2}^{l-1}  \alpha^{-}_{j,i}\right) \right) \right).\\
\end{aligned}
\end{equation*}

Using AM-GM inequality, we have
\begin{equation*}
\begin{aligned}
	&\log \left( \prod_{j=1}^{2l} \left(2\frac{\alpha(\tilde{e}_j)}{2} \prod_{i=2}^{l-1}  Q_{j,i}\right) \right) \le  \log \left( 2^{\frac{1}{l-1}} \frac{\sum_{j=1}^{2l} \left(\frac{\alpha(\tilde{e}_j)}{2} + \sum_{i=2}^{l-1}  Q_{j,i}\right)}{2l(l-1)} \right)^{2l(l-1)},
\end{aligned}
\end{equation*}
where $Q_{j,i}\in\{\alpha^{+}_{j,i},\alpha^{-}_{j,i}\}$.
Write $S_j^{+} = \frac{\alpha(\tilde{e}_j)}{2} + \sum_{i=2}^{l-1}  \alpha^{+}_{j,i}$ and $S_j^{-} = \frac{\alpha(\tilde{e}_j)}{2} + \sum_{i=2}^{l-1}  \alpha^{-}_{j,i}$ for $j=1,\ldots,2l$. Then we have the following claim.

\begin{claim}\label{claim: S_j}
	For any induced cycle $C_{2l} =\tilde{e}_1 \tilde{e}_2\ldots \tilde{e}_{2l}\tilde{e}_1 $, we have
	$\sum_{j=1}^{2l}S_j^{+} \le ml$ and $\sum_{j=1}^{2l}S_j^{-} \le ml$.
\end{claim}

Let us first assume Claim~\ref{claim: S_j} holds. Then, we have

\begin{equation*}
\begin{aligned}
	&\log (4l\Gamma) = H(e_1, \ldots, e_{2l}) \\
	\le & \log m + \sum_{C_{2l}=\tilde{e}_1 \tilde{e}_2\ldots \tilde{e}_{2l}}\frac{1}{4l\Gamma} 2 \log \left( 2^{\frac{1}{l-1}} \frac{ml}{2l(l-1)} \right)^{2l(l-1)} \\
	\le & \log \left(  2m^{l}\left( \frac{1}{2(l-1)} \right)^{l-1} \right). \\
\end{aligned}
\end{equation*}

As a conclusion, we have that $\rho(C_k, m)=\Gamma \le \left(\frac{m}{2l}\right)^{l} \left( 1+\frac{1}{l-1} \right)^{l-1} \le e\left(\frac{m}{2l}\right)^{l}= e\left(\frac{m}{k}\right)^{k/2}$ and we are done.
In the following, we give the proof of Claim~\ref{claim: S_j}.

\vspace{.2cm}
\noindent
\textbf{Proof of Claim~\ref{claim: S_j}:}

Write the induced cycle $C_{2l} = v_1v_2\ldots v_{2l}v_1$ and $\tilde{e}_i = v_{i}v_{i+1}, i=1,\ldots,2l$.
We only prove $\sum_{j=1}^{2l}S_j^{+} \le ml$ since the proof for $\sum_{j=1}^{2l}S_j^{-} \le ml$ is similar.
Note that $\alpha(\tilde{e}_j)$ and $\alpha^{+}_{j,i}~( i=2,\ldots,l-1)$ are defined by the cardinality of some edge set.
To estimate $S_j^+$, we only need to count for each edge, how does it contribute to  $S_j^+$.
In fact, we can easily verify each edge either contributes $\frac{\alpha(\tilde{e}_j)}{2} + \alpha^{+}_{j,2}$, or $\alpha^{+}_{j,i}$ to $S_j^+$ for some $3 \le i \le l-1$ since $C_{2l}$ is an induced cycle.
So a natural bound is $\sum_{j=1}^{k}S_j^+ \le 3ml$.

By a refined analysis, we can derive the claim.
We say an edge $e=uv$ is adjacent to  $w$ if either $uw \in E(G)$ or $vw \in E(G)$.
Given $e$, let $J_e = \{ j \mid \text{$v_j$ is adjacent to $e$} \}$.
Here, $e$ may intersects the cycle but does not affect our proof.



\begin{figure}[t]
	\centering         
	\begin{minipage}{0.55\linewidth}
		\centering         
		\includegraphics[width=\linewidth]{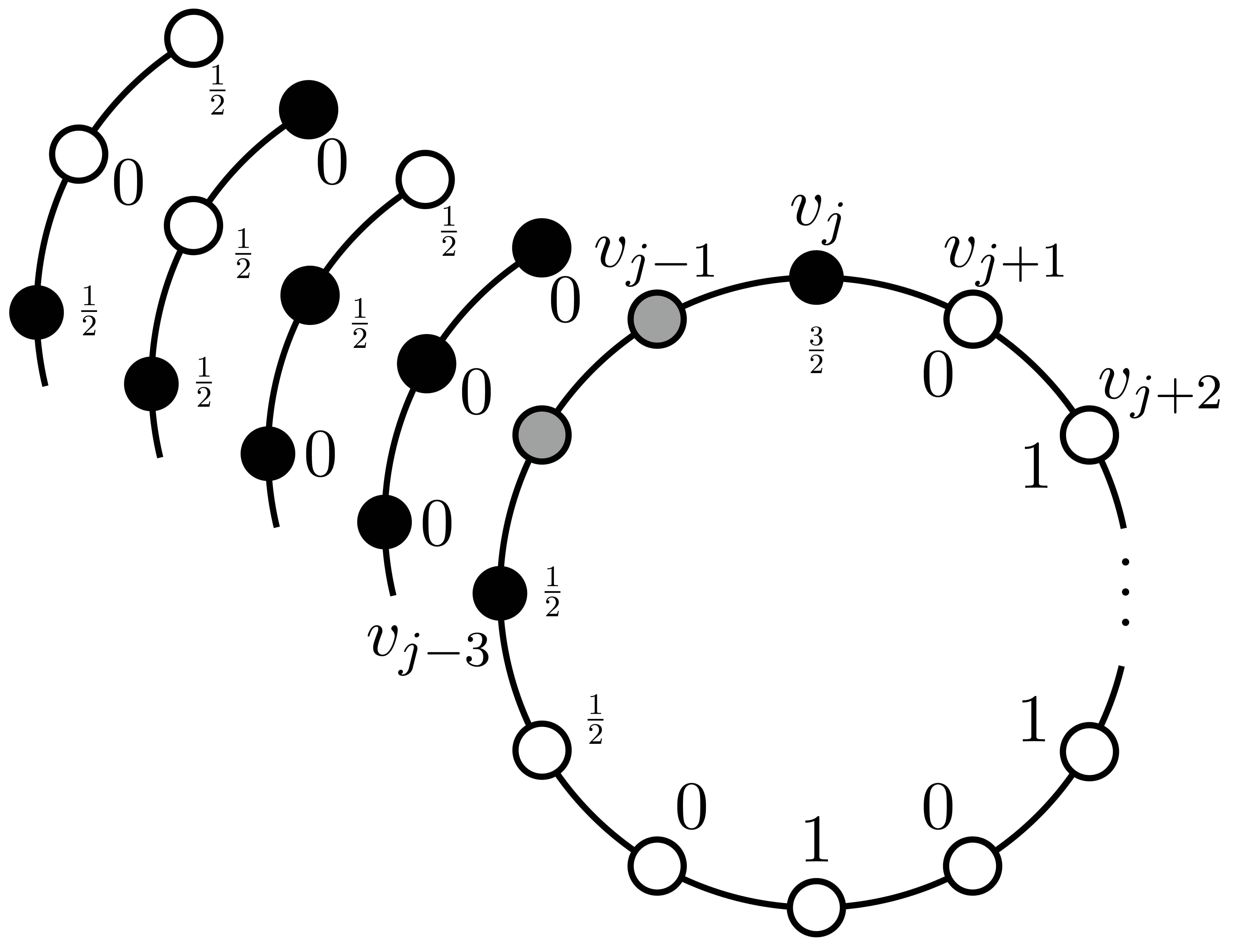}
	\end{minipage}
	\begin{minipage}{0.44\linewidth}
		\centering         
		\includegraphics[width=\linewidth]{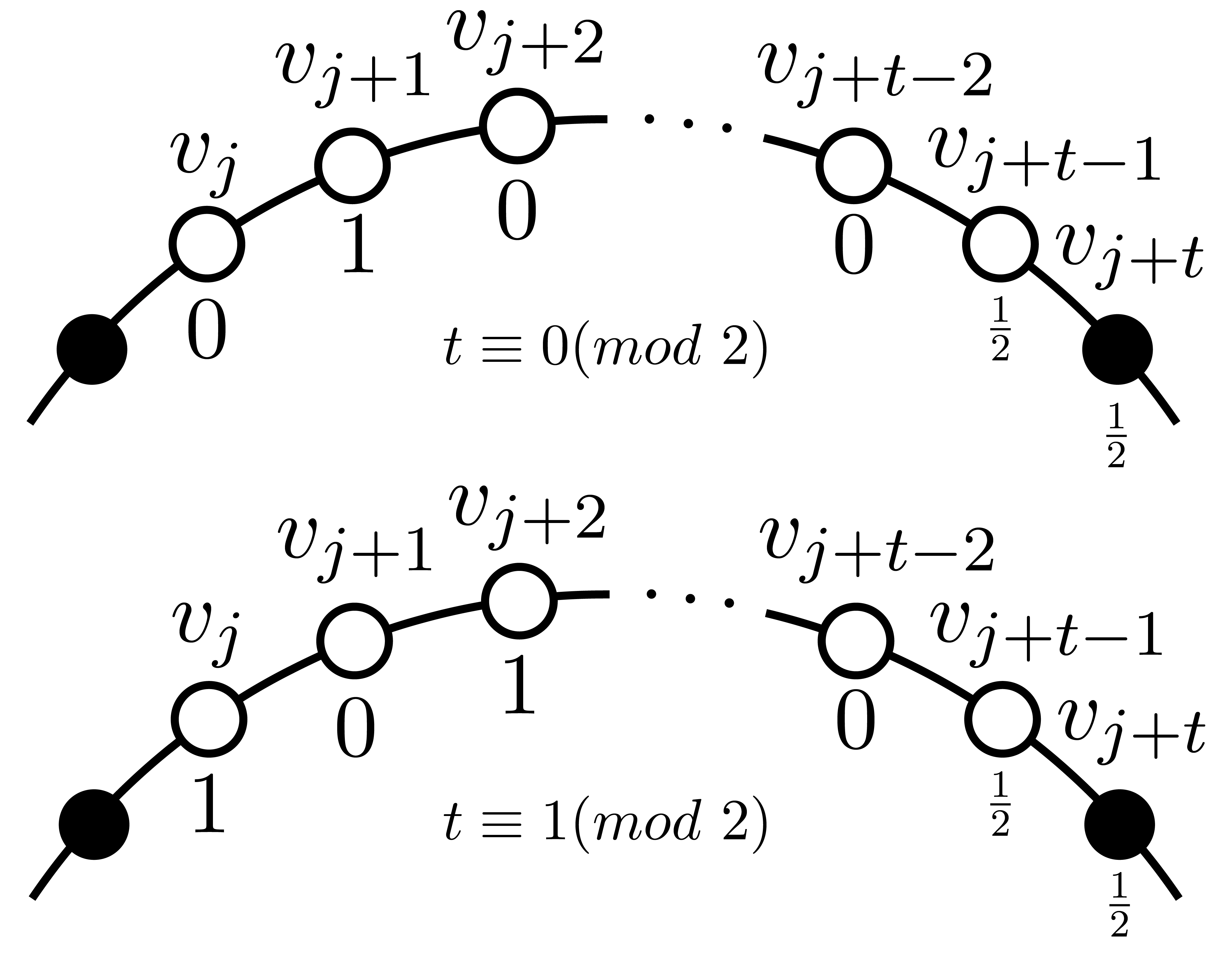}
	\end{minipage}
	\caption{The upper bound of the contribution of $e$ to $\sum_{j=1}^{2l}S_j^{+}$.
	The black nodes represent vertices in $J_e$, the white nodes represent vertices not in $J_e$ and gray nodes represent vertices which are not sure whether in $J_e$.
	The number on the side of each $v_j$ is an upper bound of the contribution of $e$ to $S_j^{+}$.
	The graph on the left shows when $e$ contributes $3/2$ to some $S_j^{+}$.
	At the upper left, we enumerate all possibilities of the gray nodes.
	The graph on the right shows the other case, split the cycle into ``white node'' paths of length $t$.   }\label{fig: Entropy}
\end{figure}

We now analyze the contribution of $e$ to $S_j^{+}$:

\begin{enumerate}
	\item When $j \in J_e$ and $j+1 \in J_e$, $e$ contributes $0$ to  $S_j^{+}$.
	\item When $j \in J_e$ and $j+1 \notin J_e$,
		\begin{enumerate}
			\item if $j+2,j+3,\ldots, j+2l-4 \notin J_e$ and $j-3 \in J_e$, $e$ contributes at most $3/2$ to $S_j^{+}$,
			\item otherwise, $e$ contributes at most $1/2$ to $S_j^{+}$.
		\end{enumerate}
	\item When $j \notin J_e$ and $j+1 \in J_e$, $e$ contributes at most $1/2$ to $S_j^{+}$.
	\item When $j,j+1 \notin J_e$, let $j'$ be the next element in $J_e$ after $j$.
		\begin{enumerate}
			\item If $|j'-j|$ is odd, $e$ contributes at most $1$ to  $S_j^{+}$.
			\item If $|j'-j|$ is even, then $e$ contributes $0$ to  $S_j^{+}$,
		\end{enumerate}
\end{enumerate}


Now consider the contribution of $e$ to $\sum_{j=1}^{2l}S_j^{+}$.
If $e$ contributes $3/2$ to $S_j^{+}$ for some $j$, then $\{j, j-3\} \subseteq J_e \subseteq \{j, j-1,j-2,j-3\}$~(see Figure~\ref{fig: Entropy}).
It is worthy to mention that if $k=6$, Claim~\ref{claim: S_j} fails on this occasion when $e$ contributes $3/2$ to $S_j^{+}$.
Enumerate all possible $J_e$ and we have $e$ contributes at most $l$ to $\sum_{j=1}^{2l}S_j^{+}$.

Otherwise, $e$ contributes at most $1$ to each $S_j^{+}$. Then we can split the cycle into maximal paths such that all vertices in the path are not in $J_e$. Depending on the parity of the length of the path, the contribution of $e$ can be estimated~(see Figure~\ref{fig: Entropy}). Sum them up then we have $e$ contributes at most $l$ to $\sum_{j=1}^{2l}S_j^{+}$.

As a conclusion, we have $\sum_{j=1}^{2l}S_j^{+}\le ml$.
\hfill$\square$

\subsection{Case 3: \texorpdfstring{$k=6$}{k=6}}

Note that Claim~\ref{claim: S_j} fails when $k=6$. If we use the method in Case 2, we can only prove $\sum_{j=1}^{k}S_j^{+} \le (l+1)m$ when $k=2l=6$, which leads to the upper bound $\rho(C_6,m) \le 4 \left( \frac{m}{6} \right)^{3}$.
Here we propose a different method to derive $\rho(C_6,m) \le 3 \left( \frac{m}{6} \right)^{3}$.

\begin{figure}[t]
	\centering         
	\includegraphics[width=0.2\linewidth]{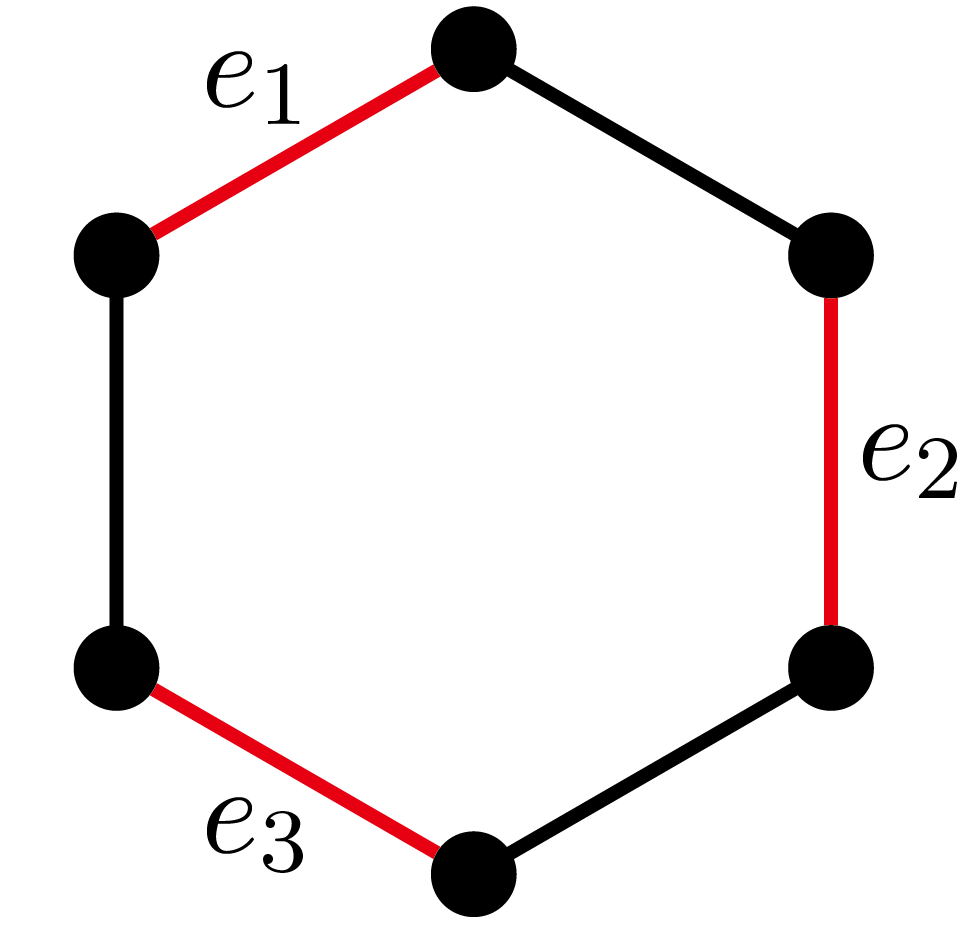}
	\caption{A sample of capable $\{e_1,e_2,e_3\}$.}\label{fig: C6 1}
\end{figure}
Recalling the definition of characterizing.
Given a graph $G$, we say $\{e_1,e_2,e_3\}$ is \textbf{capable} if $(e_1,e_2,e_3)$ characterizes an induced copy of $C_6$~(see Figure~\ref{fig: C6 1}).
Now let $G$ be the extremal graph with $m$ edges maximizing the number of induced copies of $C_{6}$ and $\Gamma$ the number of induced copies of $C_{6}$ in $G$.
Define an $3$-uniform auxiliary hypergraph $\mathcal{H}$ with vertex set $E(G)$ and edge set $\{ \{e_1,e_2,e_3\} \mid  \text{$\{e_1,e_2,e_3\}$ is capable}\}$.
Let $d(e_1,e_2)$ be the number of hyperedges in $\mathcal{H}$ containing $e_1$ and $e_2$, that is, the co-degree of $(e_1,e_2)$ in $\mathcal{H}$. Obviously, we have $e(\mathcal{H}) \triangleq |E(\mathcal{H})| = 2\Gamma$.
Let $H$ be the $2$-section of $\mathcal{H}$, that is, a graph with vertex set $V(\mathcal{H})$ and $v$ and $u$ are adjacent in $H$ if and only if there exists an hyperedge $F\in E(\mathcal{H})$ such that $v,u \in F$.
Then we have

\begin{equation*}
	\sum_{\{e_1,e_2,e_3\}\in E(\mathcal{H})}(d(e_1,e_2)+d(e_2,e_3)+d(e_1,e_3))\leq m \Gamma =\frac{1}{2}e(\mathcal{H})m,
\end{equation*}
where the inequality comes from that each edge is counted at most once in an induced copy of $C_6$. On the other hand,
\begin{equation*}
	\sum_{\{e_1,e_2,e_3\}\in E(\mathcal{H})}(d(e_1,e_2)+d(e_2,e_3)+d(e_1,e_3))=\sum_{e_1e_2\in E(H)}d(e_1,e_2)^2\geq \frac{(\sum_{e_1e_2\in E(H)}d(e_1,e_2))^2}{e(H)}.
\end{equation*}

Since $\sum_{e_1e_2\in E(H)}d(e_1,e_2)=3e(\mathcal{H}) = 6\Gamma$ and $e(H)\leq \frac{1}{2}m^2$, we have that

\[
	m\Gamma \geq \frac{{(6\Gamma)}^2}{\frac{1}{2}m^2},
\]
and then $\rho(C_6,m) = \Gamma \le  3 \left( \frac{m}{6} \right)^{3}$.





\section*{Remark}

We remark that there is a recent independent work by Chao, Antonir, Li and Yu~\cite{2025arXiv250924064C} which also study the edge inducibility problem. 
Thank them very much for giving some vital comments, including pointing out that Theorem~\ref{thm: order} is the corollary of Theorem~\ref{thm: alon friegdut kahn}.
They define the \textbf{edge inducibility} of graph $G$ by 
\[
eind(G) := \limsup_{m \to \infty} \frac{|Aut(G)|\rho(G,m)}{{(2m)}^{\alpha_f(G)}} \in [0,1].
\]
They also give the edge inducibility of $P_4$ and $P_6$ as the following theorem.
\begin{theorem}[\cite{2025arXiv250924064C}]
	\[
		eind(P_4) = \frac{1}{4}, \frac{5}{372} \le eind(P_6) \le \frac{1}{36},
	\]
\end{theorem}
which disproves our Conjecture~\ref{conj: odd path} for $k=4,6$. 
Our Conjecture~\ref{conj: cycle} coincides with their Conjecture 1.7.

\section*{Acknowledgement}

M. Lu is supported by the National Natural Science Foundation of China (Grant 12171272).

\bibliography{ref.bib}

\begin{thebibliography}{10}
\expandafter\ifx\csname urlstyle\endcsname\relax
  \providecommand{\doi}[1]{doi:\discretionary{}{}{}#1}\else
  \providecommand{\doi}{doi:\discretionary{}{}{}\begingroup \urlstyle{rm}\Url}\fi

\bibitem{alon1981number}
N.~Alon.
\newblock On the number of subgraphs of prescribed type of graphs with a given number of edges.
\newblock \emph{Israel Journal of Mathematics}, 38(1):116--130, 1981.

\bibitem{alon2016probabilistic}
N.~Alon and J.~H. Spencer.
\newblock \emph{The probabilistic method}.
\newblock John Wiley \& Sons, 2016.

\bibitem{balogh2016maximum}
J.~Balogh, P.~Hu, B.~Lidick{\`y}, and F.~Pfender.
\newblock Maximum density of induced 5-cycle is achieved by an iterated blow-up of 5-cycle.
\newblock \emph{European Journal of Combinatorics}, 52:47--58, 2016.

\bibitem{2025arXiv250924064C}
T.-W. {Chao}, A.~{Cohen Antonir}, A.~{Li}, and H.-H.~H. {Yu}.
\newblock {Edge inducibility via local directed graphs}.
\newblock \emph{arXiv e-prints}, arXiv:2509.24064, September 2025.

\bibitem{chao2024many}
T.-W. Chao, Z.~Dong, Z.~Shen, and N.~Yang.
\newblock Many cliques with small degree powers.
\newblock \emph{arXiv e-prints}, pages arXiv--2410, 2024.

\bibitem{chao2024kruskal}
T.-W. Chao and H.-H.~H. Yu.
\newblock Kruskal--katona-type problems via the entropy method.
\newblock \emph{Journal of Combinatorial Theory, Series B}, 169:480--506, 2024.

\bibitem{friedgut1998number}
E.~Friedgut and J.~Kahn.
\newblock On the number of copies of one hypergraph in another.
\newblock \emph{Israel Journal of Mathematics}, 105(1):251--256, 1998.

\bibitem{hall1987representatives}
P.~Hall.
\newblock On representatives of subsets.
\newblock \emph{Classic Papers in Combinatorics}, pages 58--62, 1987.

\bibitem{hefetz2018inducibility}
D.~Hefetz and M.~Tyomkyn.
\newblock On the inducibility of cycles.
\newblock \emph{Journal of Combinatorial Theory, Series B}, 133:243--258, 2018.

\bibitem{katona1987theorem}
G.~Katona.
\newblock A theorem of finite sets.
\newblock \emph{Classic Papers in Combinatorics}, pages 381--401, 1987.

\bibitem{kruskal1963number}
J.~B. Kruskal.
\newblock The number of simplices in a complex.
\newblock \emph{Mathematical optimization techniques}, 10:251--278, 1963.

\bibitem{KRAL2019359}
D.~Král', S.~Norin, and J.~Volec.
\newblock A bound on the inducibility of cycles.
\newblock \emph{Journal of Combinatorial Theory, Series A}, 161:359--363, 2019.
\newblock ISSN 0097-3165.
\newblock \doi{https://doi.org/10.1016/j.jcta.2018.08.003}.

\bibitem{nemhauser1975vertex}
G.~L. Nemhauser and L.~E. Trotter~Jr.
\newblock Vertex packings: structural properties and algorithms.
\newblock \emph{Mathematical Programming}, 8(1):232--248, 1975.

\bibitem{pippenger1975inducibility}
N.~Pippenger and M.~C. Golumbic.
\newblock The inducibility of graphs.
\newblock \emph{Journal of Combinatorial Theory, Series B}, 19(3):189--203, 1975.

\bibitem{willis2011bounds}
W.~Willis.
\newblock Bounds for the independence number of a graph.
\newblock 2011.

\end{thebibliography}
\bibliographystyle{wyc3}

\end{document}